\documentclass{birkmult}
\usepackage{PSFigure,psfrag}
 \theoremstyle{definition}
 
 \theoremstyle{remark}

 \numberwithin{equation}{section}

\def\IR{{\mathbb R}}
\def\IN{{\mathbb N}}

\def\IZ{{\mathbb Z}}

\def\n{\noindent}
\def\dsl{\textstyle\sum\limits}

\def\wh{\widehat}
\def\wt{\widetilde}
\def\point{{\mbox{\large $.$}}}
\def\e{\epsilon}
\def\dis{\displaystyle}

\def\f{\footnotesize}
\def\fr{\mbox{\footnotesize $\dis\frac{1}{2}$}}
\def\frvier{\mbox{\footnotesize $\dis\frac{1}{4}$}}
\def\cA{{\mathcal A}}
\def\cF{{\mathcal F}}
\def\cV{{\mathcal V}}
\def\cL{{\mathcal L}}
\def\cU{{\mathcal U}}
\def\cG{{\mathcal G}}
\def\cD{{\mathcal D}}
\def\cC{{\mathcal C}}

\renewcommand{\theequation}{\thesection.\arabic{equation}}
\newtheorem{theorem}{Theorem}[section]
\newtheorem{proposition}[theorem]{Proposition}
\newtheorem{lemma}[theorem]{Lemma}
\newtheorem{corollary}[theorem]{Corollary}
\newtheorem{remark}[theorem]{Remark}

\begin{document}
%
%
%
%
%
%
%
%
%
\title{A Lower Bound on the Disconnection Time of a Discrete Cylinder}
\author{Amir Dembo}
\address{Department of Statistics and \\
Department of Mathematics\\
Stanford University\\
Stanford, CA 94305\\
USA
}
\email{amir@math.stanford.edu}

\author{Alain-Sol Sznitman}

\address{%
Departement Mathematik\\
ETH Z\"uŸrich\\
CH 8092 Z\"urich\\
Switzerland}

\email{sznitman@math.ethz.ch}
\subjclass{60J10, 60K35, 82C41}

\keywords{disconnection time, random walk, discrete cylinders}


\begin{abstract}
We study the asymptotic behavior for large $N$ of the disconnection time $T_N$ of simple random walk on the discrete cylinder $(\IZ/N\IZ)^d \times \IZ$. When $d$ is sufficiently large, we are able to substantially improve the lower bounds on $T_N$ previously derived in \cite{DembSzni06}, for $d \ge 2$. We show here that the laws of $N^{2d}/T_N$ are tight.
\end{abstract}

\maketitle
\section{Introduction}


In this note we consider random walk on an infinite discrete cylinder having a base modelled on a $d$-dimensional discrete torus of side-length $N$. We investigate the asymptotic behavior for large $N$ of the time needed by the walk to disconnect the cylinder, or in a more picturesque language, the problem of a ``termite in a wooden beam'', see \cite{DembSzni06}, \cite{Szni06}, for recent results on this question. Building up on recent progress concerning the presence of a well-defined giant component in the vacant set left by a random walk on a large discrete torus in high dimension at times that are small multiples of the number of sites of the torus, cf.~\cite{BenjSzni06}, we are able to substantially sharpen the known lower bounds on the disconnection time when $d$ is sufficiently large.

Before describing the results of this note, we present the model more precisely. For $N \ge 1$, we consider the discrete cylinder 
\begin{equation}\label{0.1}
E = (\IZ / N \IZ)^d \times \IZ\,,
\end{equation}

\n
endowed with its natural graph structure. A finite subset $S \subseteq E$ is said to {\it disconnect} $E$ if for large $M$, $(\IZ / N \IZ)^d \times [M, \infty)$ and $(\IZ / N \IZ)^d \times (-\infty, -M]$ are contained in distinct connected components of $E \backslash S$.

We denote with $P_x$, $x \in E$, the canonical law on $E^{\IN}$ of the simple random walk on $E$ starting at $x$, and with $(X_n)_{n \ge 0}$ the canonical process. Our principal object of interest is the disconnection time:
\begin{equation}\label{0.2}
T_N = \inf\{n \ge 0; X_{[0,n]} \;\mbox{disconnects $E$}\} \,.
\end{equation}

\n
Under $P_x$, $x \in E$, the Markov chain $X_\point$ is irreducible, recurrent, and it is plain that $T_N$ is $P_x$-a.s. finite. Also $T_N$ has the same distribution under all measures $P_x$, $x \in E$. Moreover it is known, cf.~Theorem 1 of \cite{DembSzni06}, Theorem 1.2 of
\cite{Szni06} that for $d\ge1$:
\begin{equation}\label{0.3}
\lim\limits_N P_0 [N^{2d (1-\delta)} \le T_N \le N^{2d} (\log N)^{4 + \e}] = 1, \;\;\mbox{for any $\delta > 0$, $\e > 0$}\,.
\end{equation}

\n
The main object of this note is to sharpen the lower bound on $T_N$, (arguably the bound that requires a more delicate treatment in \cite{DembSzni06} and  \cite{Szni06}), when $d$ is sufficiently large. More precisely consider for $\nu \ge 1$, integer, see also (\ref{1.4}), 
\begin{equation}\label{0.4}
\mbox{$q(\nu) =$ the return probability to the origin of simple random walk on $\IZ^\nu$} \,.
\end{equation}

It is known that for large $\nu$, $q(\nu) \sim (2 \nu)^{-1}$, cf.~(5.4) in \cite{Mont56}. Our main result is:

\begin{theorem}\label{theo1.1}
Assume that $d$ is such that
\begin{equation}\label{0.5}
7 \Big(\dis\frac{2}{d+1} + \Big(1 - \dis\frac{2}{d+1}\Big) \;q(d-1)\Big) < 1\,,
\end{equation}

\n
(note that necessarily $d \ge 4$, and (\ref{0.5}) holds for large $d$, see also Remark \ref{rem2.1}), then
\begin{equation}\label{0.6}
\lim\limits_{\gamma \rightarrow 0} \;\liminf\limits_N P_0 [\gamma N^{2d} \le T_N] = 1\,.
\end{equation}

\n
(i.e. the laws of $N^{2d}/T_N$, $N \ge 2$, are tight).
\end{theorem}

In fact we expect that the laws of $\frac{T_N}{N^{2d}}$ on $(0,\infty)$ are tight, when $d \ge 2$, but the present note does not contain any novel upper bound on $T_N$ complementing (\ref{0.6}), nor a treatment of the small values of $d$.

As previously mentioned we build up on some recent progress made in the study of the vacant set left at times of order $N^{d+1}$ by a simple random walk on a $(d+1)$-dimensional discrete torus of side-length $N$, cf.~\cite{BenjSzni06}. It is shown there that when $d$ is large enough and $u$ chosen adequately small, the vacant set on the discrete torus left by the walk at time $u N^d$ contains with overwhelming probability a well-characterized giant component. In the present work we, loosely speaking, benefit for large $N$ from the presence in all blocks
\begin{equation}\label{0.7}
C_j = (\IZ/ N \IZ)^d \times \Big(\Big[\Big( j - \mbox{\f $\dis\frac{3}{4}$}\Big) \,N, \;\Big(j + \mbox{\f $\dis\frac{3}{4}$}\Big)\,N\Big] \cap \IZ\Big), \; j \in \IZ\,,
\end{equation}

\n
of such a giant component left by the walk on $E$ by time $\gamma N^{2d}$, with a probability that can be made arbitrarily close to $1$, by choosing $\gamma$ suitably small. Due to their characterization, the giant components corresponding to neighboring blocks do meet and thus offer a route left free by the walk at time $\gamma N^{2d}$, linking together ``top'' and ``bottom'' of $E$. This is in essence the argument we use when proving (\ref{0.6}).

Let us mention that the strategy we employ strongly suggests that the disconnection of $E$ takes place in a block $C_j$, which due to the presence of a sufficient number of excursions of the walk reaches criticality for the presence of a giant component. So far the study in \cite{BenjSzni06} only pertains to small values of the factor $u$ mentioned above, i.e. within the super-critical regime for the existence of a giant component in the vacant set. And the evidence of a critical and sub-critical regime, corresponding to the behavior when $u$ gets close or beyond a threshold where all components in the vacant set are typically small, rests on simulations. So the above remark is of a conjectural nature and rather points out to possible avenues of research. Incidentally progress on some of these questions may come from the investigation in \cite{Szni07} of a translation invariant model of ``random interlacements'' on $\IZ^{d+1}$, $d \ge 2$, which provides a microscopic picture of the trace left in the bulk by the walk on $E$ for the type of time regimes we are interested in.

\medskip
Let us now give more precise explanations on how we prove Theorem \ref{theo1.1}. This involves the following steps. We first analyze excursions of the walk corresponding to visits of the walk to a block of height $2N$, centered at level $j N$, thus containing $C_j$, and departures from a block with 
double height (centered at level $j N$ as well). 
We show that when $d \ge 3$, with overwhelming probability for large $N$, $C_j$ contains a wealth of segments along coordinate axes of length  $c \log N$, 
which have not been visited by the walk up to the time of  the $[u N^{d-1}]$-th excursion described above, for a small enough $u$  (see Proposition \ref{prop1.1} for a precise statement).

We then derive in Theorem \ref{theo1.2} an exponential estimate which is crucial in specifying the giant components in the various blocks $C_j$, and proving that giant components in overlapping blocks do meet. It shows that when $d \ge 4$, the probability that the walk by the time of the $[uN^{d-1}]$-th excursion has visited all points of any given set $A$ sitting in an affine-plane section of the block $C_j$, cf.~(\ref{1.12}) for the precise definition, decays exponentially in a uniform fashion with the cardinality of $A$, when $N$ is large, if $u$ is chosen small. Then a Peierl-type argument, related to the appearance of the constant $7$ in (\ref{0.5}), shows that when (\ref{0.5}) holds, with overwhelming probability the segments in $C_j$ of length $c_0 \log N$ not visited by a time $n$ prior to the $[u_0 \,N^{d-1}]$-th excursion of the walk, all belong to the same connected component of the vacant set in $C_j$ left by the walk at time $n$. Here both $c_0$ and $u_0$ are positive constants solely depending on $d$, cf.~Corollary \ref{cor1.6}.  

To prove Theorem \ref{theo1.1}, it simply remains to ascertain that with probability arbitrarily close to $1$ for large $N$, when picking $\gamma$ suitably small, no more than $u_0 N^{d-1}$ excursions have occurred in all blocks corresponding to the various levels $jN, j \in \IZ$. This is performed with the help of a coupling between the local time of simple random walk on $\IZ$ and that of Brownian motion, cf.~\cite{CsakReve83}, and a scaling argument.

\medskip
Let us now describe the structure of this note. In Section 1 we first introduce some further notation. We then show in Proposition \ref{prop1.1}, Theorem \ref{theo1.2}, Corollary \ref{cor1.6} the key ingredients for the proof of Theorem \ref{theo1.1}, i.e. the presence in a suitable regime of a multitude of segments of length $c \log N$ in the vacant set left in a block $C_j$, the above mentioned exponential bound, and the interconnection in the vacant set left in a block $C_j$ of the various segments of length $c_0 \log N$ it contains.

In Section 2 we complete the proof of Theorem \ref{theo1.1} with the help of the key steps laid out in Section 1.

Finally throughout the text $c$ or $c^\prime$ denote positive constants which solely depend on $d$, with value changing from place to place. The numbered constants $c_0,c_1,\dots,$ are fixed and refer to their first appearance in the text below. Dependence of constants on additional parameters appears in the notation. For instance $c(\gamma)$ denotes a positive constant depending on $d$ and $\gamma$.

\section{Preparation}
\setcounter{equation}{0}

In this section we first introduce some further notation and then provide with Proposition \ref{prop1.1}, Theorem \ref{theo1.2}, Corollary \ref{cor1.6}, the main ingredients for the proof of Theorem \ref{theo1.1} as explained in the Introduction.

We write $\pi_E$ for the canonical projection from $\IZ^{d+1}$ onto $E$. We denote with $(e_i)_{1 \le i \le d+1}$ the canonical basis of $\IR^{d+1}$. For $x \in \IZ^{d+1}$, resp. $x \in E$, we let $x^{d+1}$ stand for the last component, resp. the projection on $\IZ$, of $x$. We denote with $|\cdot |$ and $|\cdot |_\infty$ the Euclidean and $\ell^\infty$-distances on $\IZ^{d+1}$, or the corresponding distances induced on $E$. We write $B(x,r)$ for the closed $|\cdot|_\infty$-ball with radius $r \ge 0$, and center $x \in \IZ^{d+1}$, or $x \in E$. For $A$ and $B$ subsets of $E$ or $\IZ^{d+1}$, we write $A+B$ for the set of elements of the form $x+y$, with $x \in A$ and $y \in B$. We say that $x,y$ in $\IZ^{d+1}$ or $E$ are neighbors, resp. $\star$-neighbors if $|x-y| = 1$, resp. $|x-y|_\infty = 1$. The notions of connected or $\star$-connected subsets of $\IZ^{d+1}$ or $E$ are then defined accordingly, and so are the notions of nearest neighbor path or $\star$-nearest neighbor path on $\IZ^{d+1}$ or E. For $U$ a subset of $\IZ^{d+1}$ or $E$, we denote with $|U|$ the cardinality of $U$ and $\partial U$ the boundary of $U$:
\begin{equation}\label{1.1}
\partial U = \{ x \in U^c; \exists y \in U, \,|x-y| = 1\}\,.
\end{equation}

\n
We let $(\theta_n)_{n \ge 0}$ and $(\cF_n)_{n \ge 0}$ respectively stand for the canonical shift and filtration for the process $(X_n)_{n \ge 0}$ on $E^{\IN}$. For $U \subseteq E$, $H_U$ and $T_U$ denote the entrance time and exit time in or from $U$:
\begin{equation}\label{1.2}
H_U = \inf\{n \ge 0; \,X_n \in U\}, \;T_U = \inf\{n \ge 0; \,X_n \notin U\}\,,
\end{equation}

\n
and $\wt{H}_U$ the hitting time of $U$:
\begin{equation}\label{1.3}
\wt{H}_U = \inf\{ n \ge 1; \;X_n \in U\}\,.
\end{equation}

\n
When $U = \{x\}$, we write $H_x$ or $\wt{H}_x$ in place of $H_{\{x\}}$, or $\wt{H}_{\{x\}}$. For $\nu \ge 1$ and $x \in \IZ^\nu$, $Q^\nu_x$ denotes the canonical law on $(\IZ^\nu)^{\IN}$ of the simple random walk on $\IZ^\nu$ starting from $x$. When this causes no confusion, we use the same notation as above for the corresponding canonical process, canonical shift, the entrance, exit, or hitting times. 
So for instance, cf.~(\ref{0.4}), we have
\begin{equation}\label{1.4}
q(\nu) = Q_0^\nu [\wt{H}_0 < \infty], \;\mbox{for $\nu \ge 1$}\,.
\end{equation}

\n
We are interested in certain excursions of the walk on $E$ around the level $j N$, $j \in \IZ$, in the discrete cylinder. For this purpose we introduce for $j \in \IZ$, the blocks:
\begin{equation}\label{1.5}
\begin{split}
B_j  &= (\IZ / N\IZ)^d \times \big[(j-1)\,N, (j+1) N\big] \subseteq \wt{B}_j
\\
& = (\IZ / N \IZ)^d \times \big[(j-2) N+1, (j+2) N-1\big]\,,
\end{split}
\end{equation}

\n
so that with the definition (\ref{0.7}) we find
\begin{equation}\label{1.6}
C_j \subseteq B_j \subseteq \wt{B}_j, \;\mbox{for $j \in \IZ$}\,.
\end{equation}

\n
When $j=0$, we simply drop the subscript from the notation. We are specifically interested in the successive returns $R^j_k$, $k \ge 1$, to $B_j$, and departures $D^j_k$, $k \ge 1$, from $\wt{B}_j$:
\begin{equation}\label{1.7}
\begin{split}
R^j_1 & = H_{B_j}, \;D^j_1 = T_{\wt{B}_j} \circ \theta_{R^j_1} + R^j_1, \;\mbox{and for}\; k \ge 1\,,
\\
R^j_{k+1} & = H_{B_j} \circ \theta_{D^j_k} + D^j_k, \;D^j_{k+1} = D^j_1 \circ \theta_{D^j_k} + D^j_k \;,
\end{split}
\end{equation}
so that
\begin{equation*}
0 \le R^j_1 \le D^j_1 \le \dots \le R^j_k \le D^j_k \le \dots \le \infty\,,
\end{equation*}

\n
and for any $x \in E$, $P_x$-a.s. these inequalities are strict except maybe for the first one. We also use the convention $R^j_0 = D^j_0 = 0$, for $j \in \IZ$, as well as $R^j_t = R^j_{[t]}$, $D^j_t = D^j_{[t]}$, for $t \ge 0$.

It will be convenient in what follows to ``spread out'' the distribution of the starting point of the walk on $E$, and to this end we define:
\begin{align}
P = &\;\mbox{the law of the walk on $E$ with initial distribution,} \label{1.8}
\\
&\;\mbox{the uniform measure on $B$},\nonumber
\end{align}

\n
where we recall the convention stated below (\ref{1.6}). We write
 $E[\cdot]$ for the corresponding expectation. As noted below
 (\ref{0.2}) $T_N$ has the same distribution under any $P_x$, $x \in E$,
and it coincides with its distribution under $P$.

As we now see, when $d \ge 3$, for large $N$, with overwhelming
$P$-probability, there is in $C_j$ a wealth of segments along the
coordinate axes of length
$c \log N$ that have not been
visited by the walk up to time
 $D^j_{uN^{d-1}}$, when $u$ is small enough.
With this in mind we introduce for $K > 0$, $j \in \IZ$, $t \ge 0$, the event:
\begin{align}
\cV_{K,j,t} = \big\{& \mbox{for all $x \in  C_j$, $e \in \IZ^{d+1}$
with $|e| = 1$, for some}\label{1.9}
\\[-0.5ex]
&0 \le i < \sqrt{N}, \; X_{[0,D^j_t]} \cap \{x + (i + [0,K \log N])
\,e\big\} = \emptyset\}\,.\nonumber
\end{align}

\n
The first step on our route to the proof of Theorem \ref{theo1.1} is:
\begin{proposition}\label{prop1.1} $(d \ge 3)$

For any $K > 0$,
\begin{equation}\label{1.10}
\limsup\limits_N \;N^{-\frac{1}{4}} \;\log \sup\limits_{j \in \IZ} \;P[\cV^c_{K,j,u N^{d-1}}] < 0, \;\mbox{for small $u > 0$} \,.
\end{equation}
\end{proposition}

\begin{proof}
Theorem 1.2 of \cite{BenjSzni06} adapted to the present context states that for small $u > 0$,
$$
 \sup\limits_{j \in \IZ} \;P[\cV^c_{K,j,u N^{d-1}}] \to 0
$$
as $N \to \infty$ (take there $\beta=1/2$ and note that the dimension $d+1$ plays the role
of $d$ in \cite{BenjSzni06}). Moreover, taking $\beta_2 = \frac{1}{4}$ and $\beta_1 = \frac{1}{4} + \frac{1}{16}$ in (\ref{1.27}) of \cite{BenjSzni06}, it is not hard to verify that the proof of this theorem
actually yields the exponential decay of probabilities as in (\ref{1.10}). Indeed the probabilities in question are bounded (up to multiplicative factor $N^{d+1}$), by the sum of those in (1.49) and (1.56) of
\cite{BenjSzni06}, each of whom is shown there to be of the stated exponential decay. The
intuition behind the argument lies in a coupon-collector heuristics. Roughly speaking the strategy of the argument is the following. Given any $C_j$, $x$ in $C_j$ and coordinate direction, we consider a collection of $[N^{\beta_1}]$ segments of length $[K \log N]$ on the ``half line'' starting at $x$ with the above chosen coordinate direction, and regular interspacing  of order $[N^{\beta_1 - \beta_2}]$. We introduce a decimation process of the above collection of segments. We consider the successive excursions between times $R^j_k$ and $D^j_k$, $k \ge 1$, of the walk. At first all segments are active and we look at the first excursion visiting one of the above segments. We call it successful and take out from the list of active segments the first (active) segment, which this excusion visits. We then look for the next successful excursion visiting a segment of the list of remaining active segments. We then delete from the list of active segments the first active segment hit by this excursion. We then carry on the decimation procedure until there is no active segment left.

As in (1.49) of \cite{BenjSzni06}, one can show that when $u$ is small for large $N$, uniformly in $x \in C_j$ and in the coordinate direction, no more than $[N^{\beta_1}] - [N^{\beta_2}]$ successful excursions can occur up to time $D^j_{uN^{d-1}}$ except on a set of probability decaying exponentially in $N^{\frac{\beta_1 + \beta_2}{2}}$.

Then as in (1.56) of \cite{BenjSzni06}, one shows that during the first $[N^{\beta_1}] - [N^{\beta_2}]$ successful excursions the total number of additional segments visited after the first hit of active segments does not exceed $\frac{1}{2} \,[N^{\beta_2}]$, except on a set of probability decaying exponentially in $N^{\beta_2}$.

This enables to bound the total number of segments visited by the walk up to $D^j_{uN^{d-1}}$ by $[N^{\beta_1}] - \frac{1}{2} \,[N^{\beta_2}]$, except on a set of probability decaying exponentially in $N^{\beta_2}$. Taking into account the polynomial growth in $N$ due to the various possible choices of $x$ in $C_j$ and coordinate direction we obtain (\ref{1.10}). We refer the reader to \cite{BenjSzni06} for more details.
\end{proof}

Our next step is an exponential bound for which we need some additional notation. For $1 \le m \le d+1$, we write $\cL_m$ for the collection of subsets of $E$ that are image under the projection $\pi_E$ of affine lattices of $\IZ^{d+1}$ generated by $m$ distinct vectors of the canonical basis $(e_i)_{1 \le i \le d+1}$:
\begin{align}
\cL_m = \Big\{& \mbox{$F \subseteq E$; for some $I \subseteq \{1,\dots, d+1\}$, with $|I| = m$ and some}\label{1.11}
\\
&y \in \IZ^{d+1}, \;F = \pi_E \Big(y + \dsl_{i \in I} \;\IZ \,e_i\Big)\Big\}, \;1 \le m \le d+1\,.\nonumber
\end{align}
For $j \in \IZ$, $1 \le m \le d+1$, we consider
\begin{align}
\cA^j_m = &\; \mbox{the collection of non-empty subsets $A$ of $C_j$}\label{1.12}
\\
&\; \mbox{such that $A \subseteq F$ for some $F \in \cL_m$}\,.\nonumber
\end{align}

\n
It is plain that $\cA^j_m$ increases with $m$, and $\cA^j_{d+1}$ is the collection of non-empty subsets of 
$C_j$. 
Very much in the spirit of Theorem 2.1 of \cite{BenjSzni06}, we have the exponential bound:

\begin{theorem}\label{theo1.2} $( d \ge 3, 1 \le m \le d-2)$

Assume that $\lambda > 0$ satisfies
\begin{equation}\label{1.13}
\chi(\lambda) \stackrel{\rm def}{=} e^\lambda \;\Big(\dis\frac{m}{d+1} + \Big(1 - \dis\frac{m}{d+1}\Big) \;q(d+1 - m)\Big) < 1\,,
\end{equation}
then for $u > 0$,
\begin{equation}\label{1.14}
\limsup\limits_N \;\sup\limits_{j \in \IZ, A \in \cA^j_m} \;|A|^{-1} \log E\big[e^{\lambda \,\sum_{x \in A} \,1 \{H_x \le D^j_{uN^{d-1}}\}}\big] \le cu \;\dis\frac{e^\lambda - 1}{1 - \chi(\lambda)} \;,
\end{equation}

\n
and there exist $N_1(d,m,\lambda) > 0$, $u_1(d,m,\lambda) > 0$, such that for $N \ge N_1$:
\begin{equation}\label{1.15}
P\big[X_{[0,D^j_{u_1 N^{d-1}}]} \supseteq A\big] \le \exp \{ - \lambda |A|\}, \;\mbox{for all $j \in \IZ, \,A \in \cA^j_m$}\,.
\end{equation}
\end{theorem}

\begin{proof}
We use a variation on the ideas used in the proof of Theorem \ref{theo1.2} of \cite{BenjSzni06}. We consider for $j \in \IZ$, $A \in \cA^j_m$, $1 \le m \le d-2$, and $\lambda > 0$, the function
\begin{equation}\label{1.16}
\phi_j(z) = E_z \big[e^{\lambda \sum_{x \in A} 1\{H_x < T_{\wt{B}^j}\}}\big] \;(\ge 1), \;\mbox{for $z \in E$}\,.
\end{equation}

\n
It follows from the application of the strong Markov property at time $H_A$, that:
\begin{equation}\label{1.17}
\begin{split}
\phi_j(z) &= P_z [H_A \ge T_{\wt{B}^j}] +  
E_z [H_A < T_{\wt{B}^j}, \phi_j(X_{H_A})]
\\
& = 1 + E_z [H_A < T_{\wt{B}^j}, (\phi_j (X_{H_A}) - 1)] \,.
\end{split}
\end{equation}

\n
Now for $z \in 
\partial (B_j^c)  
\,(= (\IZ/N \IZ)^d \times \{(j - 1)N, (j+1)N\})$, we have
\begin{equation}\label{1.18}
\begin{split}
\phi_j(z) & = 1 + E_z [H_A < T_{\wt{B}^j}, \phi_j(X_{H_A}) - 1]
\\[1ex]
& \le 1 + P_z [H_A < T_{\wt{B}^j}] (\|\phi_j\|_\infty - 1)
\\ 
&\le 1 + c \;\dis\frac{|A|}{N^{d-1}} \;(\|\phi_j\|_\infty - 1) \le \exp \Big\{ c \;\dis\frac{|A|}{N^{d-1}} \;(\|\phi_j\|_\infty-1)\Big\}\,,
\end{split}
\end{equation}

\n
where in the first inequality of the last line we have used estimates on the Green function of simple random walk killed outside a strip, cf.~(2.14) of \cite{Szni03a}, to bound $P_z[H_A < T_{\wt{B}^j}]$ from above. We thus see that for $k \ge 1$,
\begin{equation}\label{1.19}
\begin{array}{l}
E\big[e^{\lambda \;\sum_{x \in A} \,1\{H_x < D^j_{k+1}\}}\big] \le E\big[e^{\lambda \sum_{x \in A} 1\{H_x < D_k^j\}}
\\[1ex]
E_{X_{R^j_{k+1}}} \big[e^{\lambda \sum_{x \in A} \,1\{H_x < T_{\wt{B}^j}\}}\big]\big] \stackrel{(\ref{1.18})}{\le} E\big[e^{\lambda \sum_{x \in A} \,1\{H_x < D^j_k\}}\big]
\\[3ex]
e^{c \frac{|A|}{N^{d-1}} \;(\|\phi_j\|_\infty - 1)}\,.
\end{array}
\end{equation}

\n
With the help of the symmetry of the Green function of the walk
killed outside $\wt{B}^j$ we show at the end of the proof
of Lemma 2.3 of \cite{DembSzni06} that
$P [H_A < D^j_1] \leq  c  |A|  N^{-(d-1)}$ in case $A$ is a sub-block
of side length $[N^\gamma]$ for fixed $0<\gamma<1$. Since exactly the
same argument (and bound) applies for all subsets $A$ of $C_j$,
we also have:
\begin{equation}\label{1.20}
\begin{array}{l}
E\big[e^{\lambda \;\sum_{x \in A} \,
1\{H_x < D^j_1\}}\big] \le 1 + P [H_A < D^j_1]
(\|\phi_j\|_\infty - 1)
\\[1ex]
\le 1 + c \;\dis\frac{|A|}{N^{d-1}} \;(\|\phi_j\|_\infty - 1)
\le e^{c \,\frac{|A|}{N^{d-1}} \; (\|\phi_j\|_\infty - 1)}\,.
\end{array}
\end{equation}

\n
Combining (\ref{1.19}) and (\ref{1.20}),
using induction as well as (\ref{1.20}) for the last term, we obtain:
\begin{equation}\label{1.21}
\begin{split}
 E\big[e^{\lambda \;\sum_{x \in A} \,1\{H_x < D^j_{uN^{d-1}}\}}\big] & \le e^{c \,\frac{|A|}{N^{d-1}} \;uN^{d-1} (\|\phi_j\|_\infty - 1)}
 \\
 & = \exp\{c u\,|A| (\|\phi_j\|_\infty - 1)\}\,.
 \end{split}
\end{equation}
We will now bound $\|\phi_j\|_\infty$.
\begin{lemma}\label{lem1.3} $(d \ge 3$, $1 \le m \le d-2$, $e^\lambda m < d+1$, $N\ge 2)$
\begin{equation}\label{1.22}
\|\phi_j\|_\infty \le \dis\frac{e^\lambda}{1 - e^\lambda \frac{m}{d+1}} \;\Big(1 - \dis\frac{m}{d+1}\Big) \;\big(1 + 
(\|\phi_j\|_\infty - 1) q_N \big) \,,
\end{equation}
where we use the notation
\begin{equation}\label{1.23}
q_N = \sup\limits_{F \in \cL_m, z \in \partial F} P_z [H_F < T_{\wt{B}^j}]\,,
\end{equation}

\n
and this quantity does not depend on $j$ due to translation invariance
of $\cL_m$ and the walk.
\end{lemma}

\begin{proof}
We consider $F \in \cL_m$, $A \subseteq F \cap (C_j \cup \partial C_j)$, and introduce the return time to $F$:
\begin{equation}\label{1.24}
R_F = H_F \circ \theta_{T_F} + T_F\,.
\end{equation}
For $z \in E$, we find:
\begin{equation}\label{1.25}
\begin{array}{l}
\phi_j(z) = E_z \big[e^{\lambda \;\sum_{x \in A} \;1\{H_x < T_{\wt{B}^j}\}}\big] \le
\\[1ex]
E_z \big[e^{\lambda (T_F + 1 \{R_F < T_{\wt{B}^j}\} 
\{(\sum_{x \in A} \,1\{H_x < T_{\wt{B}^j}\}) \circ \theta_{R_F}\} )} \big] =
\\[1ex]
E_z \big[e^{\lambda  T_F} \big(1 \{R_F \ge T_{\wt{B}^j}\}\; + 
1 \{R_F < T_{\wt{B}^j}\}\, e^{\lambda \sum_{x \in A} \, 1\{H_x < T_{\wt{B}^j}\}} \circ \theta_{R_F}\big)\big] =  
\\[2ex]
E_z \big[e^{\lambda  T_F} \big(1 + 1 \{R_F < T_{\wt{B}^j}\} \big(\phi_j (X_{R_F}) - 1\big)\big)\big] \le
\\[1ex]
E_z [e^{\lambda  T_F}] + 
E_z \big[e^{\lambda T_F}\,P_{X_{T_F}}[H_F < T_{\wt{B}^j}]\big] (\|\phi_j\|_\infty - 1) \stackrel{(\ref{1.23})}{\le}
\\[2ex]
E_z [e^{ \lambda T_F}] \,\big(1 + (\|\phi_j\|_\infty - 1) q_N \big)\,,
\end{array}
\end{equation}

\n
where we have used the strong Markov property respectively at time $R_F$ and $T_F$ in the fourth and fifth line. Then observe that when $z \notin F$, $T_F = 0$, $P_z$-a.s., whereas when $z \in F$, $T_F$ has geometric distribution with success probability $1 - \frac{m}{d+1}$, so that with $\lambda$ satisfying the hypothesis of the lemma,
\begin{equation}\label{1.26}
\begin{split}
E_z [e^{\lambda T_F}]  = &\dsl_{k \ge 1} \;\Big(1 - \dis\frac{m}{d+1}\Big) \Big(\dis\frac{m}{d+1}\Big)^{k-1} \,e^{\lambda k} =  
\\
& e^\lambda \Big(1 - \dis\frac{m}{d+1}\Big) \Big(1 - \dis\frac{e^\lambda m}{d+1}\Big)^{-1} \;.
\end{split}
\end{equation}

\n
Our claim (\ref{1.22}) follows from the last line of (\ref{1.25}).
\end{proof}

We now relate $q_N$ to $q(d+1 - m)$, cf.~(\ref{0.4}) and (\ref{1.4}). Note that our assumptions ensure that $d+1 - m \ge 3$.

\begin{lemma}\label{lem1.4} $(d \ge 3, 1 \le m \le d-2)$
\begin{equation}\label{1.27}
\limsup\limits_N \;q_N \le q(d + 1 - m)\;.\
\end{equation}
\end{lemma}

\begin{proof}
Without loss of generality we set $j =0$ in (\ref{1.23}). Then for $M \ge 1$, $F \in \cL_m$, $z \in \partial F$, we have
\begin{equation}\label{1.28}
P_z [H_F < T_{\wt{B}}] \le P_z[H_F < MN^2] + P_z [T_{\wt{B}} > MN^2]\,.
\end{equation}

\n
Using the fact that, cf.~(\ref{1.19}) of \cite{DembSzni06}:
\begin{equation}\label{1.29}
\sup\limits_{x \in \wt{B}} \;E_x \Big[\exp\Big\{\dis\frac{c}{N^2} \;T_{\wt{B}}\Big\}\Big] \le c^\prime\,,
\end{equation}

\n
to bound the last term of (\ref{1.28}), we obtain:
\begin{equation}\label{1.30}
\begin{split}
P_z [H_F < T_{\wt{B}}] & \le P_z[H_F < MN^2] + c^\prime \,e^{-cM}
\\
& \le P_{\wt{z}}^{(\IZ / N\IZ)^{d+1}} [H_{\wt{F}} < MN^2] + c^\prime \,e^{-cM}\,,
\end{split}
\end{equation}

\n
where $\wt{z}, \wt{F}$ are the respective images of $z$ and $F$ under the canonical projection from $E$ onto $(\IZ / N\IZ)^{d+1}$, and $P_{\wt{z}}^{(\IZ/N\IZ)^{d+1}}$ denotes the canonical law of simple random walk on $(\IZ / N\IZ)^{d+1}$ starting from $\wt{z}$. If we now consider the motion of the walk ``transversal to $\wt{F}$'', we find that for $N \ge2$,
\begin{equation}\label{1.31}
P_{\wt{z}}^{(\IZ/N\IZ)^{d+1}} [H_{\wt{F}} < MN^2] \le P_{e_1}^{(\IZ/N\IZ)^{d+1-m}} [H_0 < MN^2]\,,
\end{equation}

\n
with hopefully obvious notation. The right-hand side is the probability that simple random walk on $\IZ^{d+1 - m}$  starting at $e_1$ reaches $N \IZ^{d+1 - m}$ before time $MN^2$. The proof of Lemma 2.3 of \cite{BenjSzni06}  shows that the contribution of points of $N \IZ^{d+1 - m}$ other than $0$ becomes negligible as $N$ tends to infinity, so that
\begin{equation}\label{1.32}
\limsup\limits_N P_{e_1}^{(\IZ/N\IZ)^{d+1-m}}[H_0 < MN^2] \le q(d+1-m)\,,
\end{equation}
so that with (\ref{1.30}) we find
\begin{equation}\label{1.33}
\limsup\limits_N \;q_N \le q(d+1-m) + c^\prime e^{-cM}, \;\mbox{for $M \ge 1$ arbitrary}\,.
\end{equation}

\n
Letting $M$ tend to infinity, we obtain (\ref{1.27}).
\end{proof}

With (\ref{1.22}), (\ref{1.27}), it is straightforward to deduce that when $\chi(\lambda) < 1$, cf.~(\ref{1.13}),
\begin{equation}\label{1.34}
\limsup\limits_N \;\sup\limits_{j \in \IZ, A \in \cA^j_m} (\|
\phi_j \|_\infty - 1) \le \dis\frac{e^\lambda - 1}{1- \chi(\lambda)} \;.
\end{equation}

\n
Coming back to (\ref{1.21}), taking logarithms and dividing by $|A|$, the claim (\ref{1.14}) follows. As for (\ref{1.15}), we pick $\wt{\lambda}(d,m,\lambda) > \lambda$, $\wt{q} (d,m,\lambda) > q(d + 1-m)$, so that
\begin{equation}\label{1.35}
1 - e^{\wt{\lambda}} \Big(\dis\frac{m}{d+1} + \Big(1 - \dis\frac{m}{d+1}\Big) \,\wt{q}\Big) = \fr \;\big(1 - \chi(\lambda)\big)\,.
\end{equation}

\medskip\n
Applying (\ref{1.14}) with $\wt{\lambda}$ (which satisfies $\chi(\wt{\lambda}) < 1)$, we see that for $u > 0$, $N \ge N_2(d,m,\lambda,u)$, and any $j \in \IZ$, $A \in \cA^j_m$, one has with (\ref{1.35}):
\begin{equation}\label{1.36}
\begin{split}
P[X_{[0,D^j_{uN^{d-1}}]} \supseteq A] & \le P\Big[\dsl_{x \in A} \;1 \{H_x < H_{D^j_{uN^{d-1}}}\} \ge |A| \Big]
\\
& \le \exp\Big\{- \wt{\lambda} \,|A| + c u \;\dis\frac{e^{\wt{\lambda}} - 1}{1 - 
 \chi(\wt\lambda)} \;|A|\Big\}\,.
\end{split}
\end{equation}

\n
Choosing $u = u_1(d,m,\lambda)$ small enough, and setting $N_1(d,m,\lambda) = N_2(d,m,\lambda,u_1)$, we obtain (\ref{1.15}).
\end{proof}

We will now use the above exponential control combined with a Peierl-type argument to ensure the typical presence for large $N$ of a well-specified giant component in the vacant set left by the walk in a block $C_j$, as long as the number of excursions between $B_j$ and $\wt{B}_j$ does not exceed a small multiple of $N^{d-1}$. This construction will force giant components corresponding to neighboring blocks to have non-empty intersection. We recall that $\star$-nearest neighbor paths have been defined at the beginning of this section, and introduce
\begin{align}
a(n) = &\; \mbox{the cardinality of the collection of $\star$-nearest neighbor} \label{1.37}
\\[-0.5ex]
&\;\mbox{self-avoiding paths on $\IZ^2$, starting at the origin, with $n$ steps}\,.\nonumber
\end{align}

\n
One has the straightforward upper bound
\begin{equation}\label{1.38}
a(n) \le 8  \; 7^{n-1}, \;\mbox{for $n \ge 1$}\,.
\end{equation}

\n
Given $N \ge 1$, $K > 0$, $j \in \IZ$, $t \ge 0$, we introduce the event
\begin{align}
\cU_{K,j,t} = &\; \mbox{for any $F \in \cL_2$, $n \le D^j_t$, any connected subsets $O_1,O_2$} \label{1.39}
\\
&\;\mbox{of $F \cap C_j \backslash X_{[0,n]}$, with $|\cdot |_\infty$-diameter at least $[K \log N]$}\nonumber
\\
&\;\mbox{are in the same connected component of $F \cap C_j \backslash X_{[0,n]}$}\,.\nonumber
\end{align}

\n
This event will be helpful in specifying the above mentioned giant components. We recall the notation (\ref{1.4}).

\begin{corollary}\label{cor1.4}
If $d \ge 4$ is such that
\begin{equation}\label{1.40}
\rho \stackrel{\rm def}{=} 7 \;\Big(\dis\frac{2}{d+1} + \Big(1 - \dis\frac{2}{d+1}\Big) \;q(d-1)\Big) < 1\,,
\end{equation}

\n
as this happens for any large $d$, then there are constants $c_0 > 0$, cf.~(\ref{1.45}), and $u_0$, cf.~(\ref{1.43}), such that
\begin{equation}\label{1.41}
\limsup\limits_N \;N^{2d} \;\sup\limits_j \;P [\cU^c_{c_{0,j,u_o N^{d-1}}}] = 0\,.
\end{equation}
\end{corollary}

\begin{proof}
Note that $q(\nu) \sim (2 \nu)^{-1}$, cf.~(5.4) of \cite{Mont56}, and (\ref{1.40}) holds for any large enough $d$. Assume that (\ref{1.40}) holds and choose $\lambda_0(d)$ such that
\begin{equation}\label{1.42}
e^{\lambda_0} = 7 \rho^{-\frac{1}{2}} > 7, \;\mbox{so that}\;e^{\lambda_0} \Big(\dis\frac{2}{d+1} + \Big(1 - \dis\frac{2}{d+1}\Big) \;q(d-1)\Big) = \rho^{\frac{1}{2}} < 1\,.
\end{equation}

\n
When $N$ is large, on $\cU^c_{K,j,t}$, one can find $F \in \cL_2$, $n \le D^j_t$, $O_1,O_2 \subseteq F \cap C_j \backslash X_{[0,n]}$, distinct connected components of $F \cap C_j \backslash X_{[0,n]}$ with $|\cdot|_\infty$-diameter at least $[K \log N]$. If the last vector $e_{d+1}$ of the canonical basis does not enter the definition of $F$, cf.~(\ref{1.11}), then $F \subseteq C_j$, and we can introduce an affine projection of $\IZ^2$ onto $F$, and define $\wh{O}_i$, $i=1,2$, the inverse images of $O_i$ under this affine projection. Considering separately the case when at least one of the $\wh{O}_i$, $i=1,2$, has bounded components, (necessarily of $|\cdot |_\infty$-diameter at least $[K \log N]$), or both of the $\wh{O}_i$ have unbounded components, one can construct a $\star$-nearest neighbor self-avoiding path $\pi$ with $[K \log N]$ steps in $\partial O_1 \cap F$ or $\partial O_2 \cap F  \subseteq F \cap X_{[0,n]} \subseteq F \cap C_j \cap X_{[0,D^j_t]}$, see also Proposition 2.1, p.~387, in \cite{Kest82}. On the other hand if the last vector $e_{d+1}$ of the canonical basis enters the definition of $F$, we introduce an affine projection of $\IZ^2$ onto $F$ so that the inverse image of $F \cap C_j$ coincides with the strip $\IZ \times ([-\frac{3}{4} N, \frac{3}{4} N] \cap \IZ)$. Defining as above $\wh{O}_i$, $i=1,2$, the inverse images of $O_i$ under this affine projection, we can separately consider the case when at least one of the $\wh{O}_i$, $i=1,2$ has bounded components, (necessarily of $|\cdot|_\infty$-diameter at least $[K \log N]$), or both of the $\wh{O}_i$ have unbounded components. We can  then construct a $\star$-nearest neighbor self-according path $\pi$ with $[K \log N]$ steps in $\partial O_1 \cap F \cap C_j$ or $\partial O_2 \cap F$  $\cap C_j \subseteq$ $F \cap C_j \cap X_{[0,n]} \subseteq F \cap C_j \cap X_{[0,D^j_t]}$.

\begin{center}
\PS{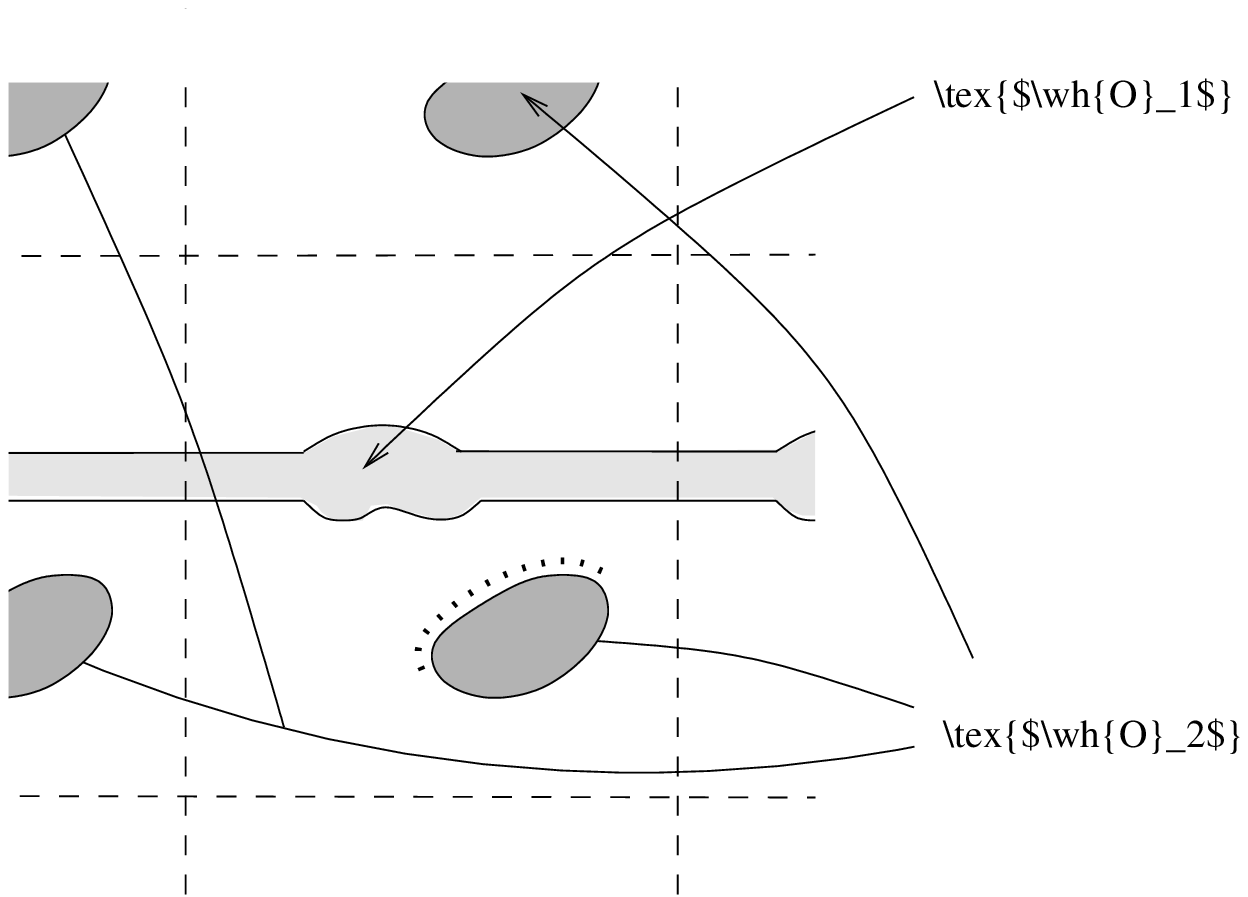}{800}
\end{center}

{\small
\begin{center}
\begin{tabular}{lp{12cm}}
Fig. 1:& An example of possible $\wh{O}_1$, $\wh{O}_2$ is depicted in the case when $e_{d+1}$ does not enter the definition of $F$. The square deliminited by dashed lines is a ``fundamental domain'' for the affine projection. The dotted line is $*$-connected and projects on a subset of $F \cap C_j \cap X_{[0,n]}$.
\end{tabular}
\end{center}
}

As a result setting, cf.~above (1.15),
\begin{equation}\label{1.43}
u_0(d)  = u_1\big(d,m=2,\lambda = \lambda_0(d)\big)\,,
\end{equation}
we see that for large enough $N$, for any $j \in \IZ$,
\begin{equation}\label{1.44}
\begin{array}{lcl}
P[\cU^c_{K,j,u_0 N^{d-1}}] &\hspace{-1ex} \le&\hspace{-1ex} \dsl_F \; \dsl_\pi \;P\big[X_{[0,D^j_{u_0 N^{d-1}}]}\supseteq A\big]
\\[1ex]
&\hspace{-1ex} \stackrel{(\ref{1.15})}{\le}&\hspace{-1ex} \sup\limits_j \;\dsl_F \;\dsl_\pi \;e^{-\lambda_0|A|}
\\[1ex]
&\hspace{-1ex} \stackrel{(\ref{1.38})}{\le}&\hspace{-1ex} \sup\limits_j \;\dsl_F \;c N^2 \, 7^{[K \log N]-1} \,e^{-\lambda_0[K \log N]}
\\[1ex]
&\hspace{-1ex} \stackrel{(\ref{1.42})}{\le}&\hspace{-1ex} c\,N^{d+1} \rho^{\frac{1}{2} [K \log N]} \;,
\end{array}
\end{equation}

\n
where the sum over $F$ pertains to $F \in \cL_2$ with $F \cap C_j \not= \phi$, the sum over $\pi$ pertains to the collection of $\star$-nearest neighbor self-avoiding paths with values in $F \cap C_j$ with $[K \log N]$ steps, and $A$ stands for the set of points visited by $\pi$. If we now specify $K$ to take the value
\begin{equation}\label{1.45}
c_0 = 8d \Big( \log \dis\frac{1}{\rho}\Big)^{-1}\,,
\end{equation}

\n
the claim (\ref{1.41}) follows from the last line of (\ref{1.44}).
\end{proof}

We now introduce for $j \in \IZ, t \ge 0$, the event, cf.~(\ref{1.9}), (\ref{1.39}):
\begin{equation}\label{1.46}
\cG_{j,t} = \cV_{c_{0,j,t}} \cap \cU_{c_{0,j,t}} \,.
\end{equation}

\begin{corollary}\label{cor1.6}
Assume $d \ge 3$, and $N \ge 2$ large enough so that $(\IZ/N\IZ)^d$ has $|\cdot|_\infty$-diameter bigger than $c_0 \log N$. Then for $j \in \IZ$, $t \ge 0$, on $\cG_{j,t}$, for any $0 \le n \le D^j_t$,
\begin{align}
&\mbox{all segments in $C_j \backslash X_{[0,n]}$ of length $L_0 \stackrel{\rm def}{=} [c_0 \log N]$ belong to}\label{1.47}
\\
&\mbox{the same connected component of $C_j \backslash X_{[0,n]}$},\nonumber
\\[2ex]
&\mbox{any $F \in \cL_1$ intersecting $C_j$ contains a segment of length $L_0$} \label{1.48}
\\
&\mbox{included in $C_j \backslash X_{[0,n]}$}\,.\nonumber
\end{align}
Moreover when $d \ge 4$ satisfies (\ref{1.40}), with the notation (\ref{1.43}), one has:
\begin{equation}\label{1.49}
\lim\limits_N \;N^{2d} \;\sup\limits_{j \in \IZ} P[\cG^c_{j,u_0 N^{d-1}}] = 0 \,.
\end{equation}
\end{corollary}

\begin{proof}
The claim (\ref{1.48}) readily follows from the inclusion $\cG_{j,t} \subseteq \cV_{c_0,j,t}$, and (\ref{1.9}). To prove (\ref{1.47}), consider $n \le D^j_t$ and note that $\cG_{j,t} \subseteq \cU_{c_0,j,t}$. Hence any two segments of length $L_0$ contained in $\wt{F} \cap C_j \backslash X_{[0,n]}$, with $\wt{F} \in \cL_2$, belong to the same connected component of $F \cap C_j \backslash X_{[0,n]}$. Then consider $\wt{F}, \wt{F}_2 \in \cL_2$, we see that with (\ref{1.48}) and the above,
\begin{align}
&\mbox{when $\wt{F}_1 \cap \wt{F}_2 \in \cL_1$ and intersects $C_j$, all segments of length $L_0$ in}\label{1.50}
\\
&\mbox{$(\wt{F}_1 \cup \wt{F}_2) \cap C_j \backslash X_{[0,n]}$ are in the the same connected component}\nonumber
\\
&\mbox{of $C_j \backslash X_{[0,n]}$}\,.\nonumber
\end{align}

\n
Next, given $y_0$ and $y$ in $C_j$, we can construct a nearest neighbor path $(y_i)_{0 
\le i \le m}$ in $C_j$, with $y_m = u$. Given $\wt{F} \ni y_0$, with $\wt{F} \in \cL_2$, we can construct a sequence $\wt{F}_i \in \cL_2$, $0 \le i \le m$, such that
\begin{align}
&\mbox{$\wt{F}_0 = \wt{F}$, $y_i \in \wt{F}_i$, for $0 \le i \le m$, and either $\wt{F}_{i-1} = \wt{F}_i$, or} \label{1.51}
\\
&\mbox{$\wt{F}_{i - 1} \cap \wt{F}_i \in \cL_1$ and intersects $C_j$, for $1 \le i \le m$}, \nonumber
\end{align}

\n
(see for instance below (2.60) of \cite{BenjSzni06} for a similar argument).

Analogously when $\wt{F}, \wt{F}^\prime \in \cL_2$ have a common point $y \in C_j$, we can define $\wt{F}_i \in \cL_2$, $0 \le i \le 2$, such that
\begin{align}
&\mbox{$\wt{F}_0 = \wt{F}$, $\wt{F}_2 = \wt{F}^\prime$, with $y \in \wt{F}_i$, $0 \le i \le 2$, and either $\wt{F}_i = \wt{F}_{i-1}$ or} \label{1.52}
\\
&\mbox{$\wt{F}_i \cap \wt{F}_{i-1} \in \cL_1$ (and intersects $C_j$), for $i = 1,2$}.\nonumber
\end{align}

\n
Combining (\ref{1.48}) and (\ref{1.50}) - (\ref{1.52}), we obtain (\ref{1.47}). Finally 
(\ref{1.49}) is a direct consequence of (\ref{1.10}) and (\ref{1.41}).
\end{proof}

\section{Denouement}
\setcounter{equation}{0}

We now use the results of the previous section to prove Theorem \ref{theo1.1}. As mentioned in the Introduction, the rough idea is that by making $\gamma$ small one can ensure that with arbitrarily high probability, when $N$ is large enough, for all $j$ in $\IZ$ the times $D^j_{u_0 N^{d-1}}$ are bigger than $\gamma N^{2d}$. Then on ``most'' of the event $\{\inf_{j \in \IZ} D^j_{u_0 N^{d-1}} > \gamma N^{2d}\}$, there is a profusion of segments of length $L_0$ in each $C_j \backslash X_{[0,[\gamma N^{2d}]]}$, and they all lie in the same connected component of $X^c_{[0,[\gamma N^{2d}]]}$. This now forces the disconnection time $T_N$ to be bigger than $\gamma N^{2d}$.

\bigskip\n
{\it Proof of Theorem \ref{theo1.1}.} As note below (\ref{0.2}), we can with no loss of generality replace $P_0$ with $P$ in the claim (\ref{0.6}) to be proved. We introduce for $0 < \gamma < 1$, $t \ge 0$, the events:
\begin{align}
\cC_{\gamma,t} & = \cD_{\gamma,t} \cap \big(\bigcap_{|j| \le 2N^{2d-1}} \cG_{j,t}\big), \;\mbox{where} \label{2.1}
\\[1ex]
\cD_{\gamma,t} & =\bigcap_{|j| \le 2N^{2d-1}} \; \{D^j_t > \gamma N^{2d}\}\,. \label{2.2}
\end{align}

\n
Note that $P$-a.s., the vertical component of the walk up to time $N^{2d}$ remains in $[-N^{2d} - N, N^{2d} + N]$. Hence for large $N$, with (\ref{1.47}),  (\ref{1.48}), $P$-a.s. on $\cC_{\gamma_,t}$, there is a nearest neighbor path in  $X^c_{[0,[\gamma N^{2d}]]}$ starting in $(\IZ / N \IZ)^d \times (-\infty, - M]$ and ending in $(\IZ / N \IZ)^d \times [M, + \infty)$ for any $M \ge 1$, and therefore $T_N > \gamma N^{2d}$. As a result Theorem \ref{theo1.1} will be proved once we show that with the notation of (\ref{1.43})
\begin{align}
&\liminf\limits_N P\big[ \bigcap_{|j| \le 2N^{2d-1}}\, \cG_{j,u_0 N^{d-1}}\big] = 1, \;\mbox{and} \label{2.3}
\\[1ex]
&\lim\limits_{\gamma \rightarrow 0}\; \liminf\limits_N P [\cD_{\gamma, u N^{d-1}}]  = 1, \;\mbox{for any $u > 0$}\,. \label{2.4}
\end{align}

\n
The claim (\ref{2.3}) readily follows from (\ref{1.49}), so we only need to prove (\ref{2.4}). To this end we introduce the following sequence of $(\cF_n)$-stopping times that are $P_x$-a.s. finite for any $x \in E$:
\begin{equation}\label{2.5}
\begin{split}
\tau_0 & = H_{(\IZ / N \IZ)^d \times  N \IZ}, \; \mbox{cf.~(\ref{1.2}) for the notation, and for $k \ge 0$,}
\\[1ex]
\tau_{k+1} & = \inf \{n > \tau_k; \,|X_n^{d+1} - X^{d+1}_{\tau_k} | = N\}\,,
\end{split}
\end{equation}

\n
where according to the notation of the beginning of Section 1, $X_\point^{d+1}$ denotes the $\IZ$-component of $X_\point$. We also write $\tau_t = \tau_{[t]}$, for $t \ge 0$. The count of visits of $X_{\tau_n}$, $n \ge 0$, to level $\ell N$ with $n$ at most $t$, is then expressed by
\begin{equation}\label{2.6}
L_N (\ell,t) = \dsl_{0 \le n \le t} 1\{X^{d+1}_{\tau_n} = \ell N\}, \;\mbox{for $\ell \in \IZ, t \ge 0$}\,.
\end{equation}

\n
Pick $M \ge 1$, and note that as soon as $\tau_{M \gamma N^{2d-2}} > \gamma N^{2d}$ and $L_N (\ell, M \gamma N^{2d-2}) < \frac{1}{2} \,[uN^{d-1}]$, for all $|\ell | \le 2 N^{2d-1} + 1$, then $D^j_{uN^{d-1}} > \gamma N^{2d}$ for all $|j| \le 2N^{2d-1}$, and hence $\cD_{\gamma, uN^{d-1}}$ occurs. As a result we find that
\begin{equation}\label{2.7}
\begin{array}{l}
P[\cD_{\gamma,uN^{d-1}}] \ge 
\\[1ex]
P\Big[\theta^{-1}_{\tau_0} \Big(\Big\{ \sup\limits_{|\ell | \le 2N^{2d-1} + 1} \,L_N (\ell, M \gamma N^{2d-2}) < \fr \,[uN^{d-1}]\Big\}\Big) \cap  
\\
\qquad \qquad \big\{\tau_{M \gamma N^{2d - 2}} \circ \theta_{\tau_0} > \gamma N^{2d}\big\}\Big] \ge
\\[1ex]
P_0 \Big[\Big\{\sup\limits_{|\ell | \le 2N^{2d-1} + 2} \, L_N (\ell, M \gamma N^{2d-2}) < \fr \,[uN^{d-1}]\Big\} \cap 
\\[1ex]
\qquad \qquad \{\tau_{M \gamma N^{2d - 2}}  > \gamma N^{2d}\}\Big] \ge  a_1 - a_2, \;\mbox{with the notation}
\end{array}
\end{equation}
\begin{equation}\label{2.8}
\begin{array}{l}a_1 = P_0 \Big[\sup\limits_{|\ell | \le 2N^{2d-1} + 2}  \,L_N (\ell, M \gamma N^{2d-2}) < \fr \,[uN^{d-1}]\Big]\,,
\\
\\[-2ex]
a_2 = P_0\big[\tau_{M \gamma N^{2d-2}} \le \gamma N^{2d}\big]\,,
\end{array}
\end{equation}

\medskip\n
and where we have used the strong Markov property at time $\tau_0$, as well as translation invariance when  going from the second to the third line of (\ref{2.7}).

Note that under $P_0$, $\tau_k$, $k \ge 0$, has stationary independent increments and applying the invariance principle to the $\IZ$-component of $X$, we also have
\begin{equation}\label{2.9}
E_0\Big[\exp\Big\{ - \dis\frac{c_1}{N^2} \;\tau_1\Big\}\Big] \le e^{-c_2} \,.
\end{equation}

\n
With Cheybyshev's inequality we thus find
\begin{equation*}
\begin{split}
P_0[\tau_{M \gamma N^{2d-2}} \le \gamma N^{2d}] & \le e^{c_1 \gamma N^{2d-2}} E_0\big[e^{- \frac{c_1}{N^2} \;\tau_{M \gamma N^{2d-2}}}\big]
\\[1ex]
& \le c \,e^{(c_1 - c_2 M) \gamma N^{2d-2}} \,.
\end{split}
\end{equation*}

\n
Choosing from now on $M > \frac{c_1}{c_2}$, we find that 
\begin{equation}\label{2.10}
\lim\limits_N \,a_2 = 0\,.
\end{equation}

\n
Coming back to $a_1$, we observe that under $P_0$, $L_N (\cdot,\cdot)$ has the same distribution as the local time process of simple random walk on $\IZ$ starting at the origin.

In fact, cf.~(1.20) of \cite{CsakReve83}, we can construct on some auxiliary probability space $(\wt{\Sigma}, \wt{\cA}, \wt{P})$ a one-dimensional Brownian motion $(\wt{B}_t)_{t \ge 0}$, and a simple random walk on $\IZ$ starting at the origin, $Z_k, k \ge 0$, so that setting $\wt{L}(x,t)$, $x \in \IR$, $t \ge 0$,  be a jointly continuous version of the local time of $\wt{B}_\point$ and
\begin{equation}\label{2.11}
L(x,k) = \dsl^k_{n=0} \;1\{Z_n = x\}\,, x \in \IZ, \, k \ge 0\,,
\end{equation}

\n
be the local time of the simple random walk $Z_\point$, one has
\begin{equation}\label{2.12}
\mbox{$\wt{P}$-a.s., for all $\rho > 0, \lim\limits_{n \rightarrow \infty} \,n^{-\frac{1}{4} - \rho} \,\sup\limits_{x \in \IZ} \,|\wt{L} (x,n) - L(x,n)| = 0$}\,.
\end{equation}

With this we find that for any $\gamma > 0, u > 0$,
\begin{equation*}
\begin{array}{lcl}
\liminf\limits_N \,a_1 &\hspace{-3ex} \ge &\hspace{-3ex} \liminf\limits_N \,\wt{P} \Big[ \sup\limits_{|\ell | \le 3 N^{2d-1}} \,L(\ell, [M \gamma N^{2d-2}]) < \fr \;[uN^{d-1}]\Big]
\\[1ex]
&\hspace{-3ex} \stackrel{(\ref{2.12})}{\ge} &\hspace{-2ex} \liminf\limits_N \,\wt{P} \,\Big[\sup\limits_{w \in \IR} \,\wt{L} (w,M \gamma N^{2d-2}) < \frvier \;u N^{d-1}\Big]
\\[2ex]
& \hspace{-3ex}\stackrel{\rm scaling}{=} & \hspace{-1ex}\wt{P} \Big[ \sup\limits_{v \in \IR} \,\wt{L} (v,M \gamma) < \frvier \;u \Big] \,.
\end{array}
\end{equation*}

\n
It thus follows from the $\wt{P}$-a.s. joint continuity of $\wt{L}(\cdot,\cdot)$ that
\begin{equation}\label{2.13}
\lim\limits_{\gamma \rightarrow 0} \;\liminf\limits_N \,a_1 = 1\,.
\end{equation}

\n
Together with (\ref{2.7}) and (\ref{2.10}), this concludes the proof of Theorem \ref{theo1.1}. \hfill $\square$

\begin{remark}\label{rem2.1} \rm ~
It is plain that the condition (\ref{0.5}) under which 
Theorem \ref{theo1.1} holds requires $d  \ge 14$. 
Using a numerical evaluation of $q(\nu)$ 
in (\ref{0.4}) based on the formula (5.1) of \cite{Mont56}, 
which was kindly provided to us by 
Wesley P. Petersen, one can see that (\ref{0.5}) holds 
when $d  \ge 17$ and fails when $d  < 17$. 
On the other hand the numerical simulations performed in the context 
of the investigation of the vacant set left by random walk on the 
discrete torus in \cite{BenjSzni06} make it plausible that the 
conclusion of Theorem \ref{theo1.1} should hold for all $d  \ge 1$ 
(the case $d=1$ being trivial).
\hfill $\square$

\end{remark}


\pagebreak
\subsection*{Acknowledgments}
We wish to thank Wesley P. Petersen for kindly communicating to us a table of numerical values of the return probability of simple random walk to the origin.

Amir Dembo would like to thank the FIM  for hospitality and financial support during his visit to  ETH. His research was also partially supported by the NSF grants DMS-0406042, DMS-FRG-0244323.

\end{document}